\documentclass[10pt,twoside,english]{amsart2000}
\advance\oddsidemargin by -1.0cm
\advance\evensidemargin by -1.0cm
\advance\topmargin by -1.0cm 
\usepackage{amssymb}
\usepackage{babel}
\usepackage{amstext2000}
\usepackage{amscd2000}   
\usepackage{epsfig}  
\usepackage{rotating}
\theoremstyle{plain}
\newtheorem{co}{Corollary}[section]

\newtheorem{thm}{Theorem}[section]
\newtheorem{pro}{Proposition}[section]
\theoremstyle{definition}
\newtheorem{exa}{Example}[section]

\def\bea{\begin{eqnarray}}
\def\eea{\end{eqnarray}}
\def\beann{\begin{eqnarray*}}
\def\eeann{\end{eqnarray*}}
\def\beq{\begin{equation}}
\def\eeq{\end{equation}}
\def\ba{\begin{array}}
\def\ea{\end{array}}
\def\ben{\begin{enumerate}}
\def\een{\end{enumerate}}

\date{\today}
\begin{document}
%---------------------------------------------------------------------------
\def\haken{\mathbin{\hbox to 6pt{%
                 \vrule height0.4pt width5pt depth0pt
                 \kern-.4pt
                 \vrule height6pt width0.4pt depth0pt\hss}}}
    \let \hook\intprod
%
%----------------------------------------------------------------------------
\title{Killing spinor equations in dimension $7$ and geometry of 
integrable $\mbox{G}_2$-manifolds}
%----------------------------------------------------------------------------
%
% author and address
%
%-------------------------------------------
%
\author{Thomas Friedrich and Stefan Ivanov}
%-------------------------------------------
\address{\hspace{-5mm} 
Thomas Friedrich\newline
Institut f\"ur Reine Mathematik \newline
Humboldt-Universit\"at zu Berlin\newline
Sitz: WBC Adlershof\newline
D-10099 Berlin, Germany\newline
{\normalfont\ttfamily friedric@mathematik.hu-berlin.de}}
\address{\hspace{-5mm} 
Stefan Ivanov\newline
Faculty of Mathematics and Informatics\newline
University of Sofia ``St. Kl. Ohridski''\newline
blvd. James Bourchier 5\newline
1164 Sofia, Bulgaria\newline
{\normalfont\ttfamily ivanovsp@fmi.uni-sofia.bg}}
%------------------------------------------------------
\thanks{Supported by the SFB 288 "Differential geometry
and quantum physics" of the DFG and the European Human Potential Program EDGE, 
Research Training Network HPRN-CT-2000-00101. S.Ivanov thanks ICTP for the
support and excellent environment.}
%------------------------------------------------------
%----------------------------------------------
\begin{abstract} We compute the scalar curvature of $7$-dimensional $\mbox{G}_2$-manifolds 
admitting a connection with totally skew-symmetric torsion. We prove the
formula for the general solution of the Killing spinor equation and
express the Riemannian scalar curvature of the solution in terms of the
dilation function and the NS $3$-form field. In dimension $n=7$ the dilation 
function involved
in the second fermionic string equation has an interpretation as a 
conformal change of the underlying integrable $\mbox{G}_2$-structure into a cocalibrated
one of pure type $W_3$. 
%---------------
\end{abstract}
%-------------
\maketitle
%----------------
\tableofcontents
%----------------
\pagestyle{headings}
%
%-----------------------------------------------------------------------------
\section{Introduction}\noindent
%-----------------------------------------------------------------------------
%
Riemannian manifolds admitting parallel spinors with respect to a metric connection with 
totally skew-symmetric torsion became a subject of interest in theoretical and 
mathematical physics recently. One of the main reasons is that the number of preserved 
supersymmetries in string theory depends essentially on the number of parallel spinors. 
In $10$-dimensional string theory, the Killing spinor equations with 
non-constant dilation $\Phi$ and the 3-form field strength $H$
can be written in the following way \cite{Stro}, (see  \cite{IP1,IP,FI})
\beq\label{ks}\tag{$*$}
\nabla \Psi\, =\, 0,\quad
(d\Phi - \frac{1}{2}H)\cdot\Psi\, =\, 0,
\eeq
where $\Psi$ is a spinor field and $\nabla$ is a metric connection with 
totally skew-symmetric torsion $T=H$.
The existence of a  parallel spinor imposes restrictions on the holonomy 
group since the spinor holonomy representation has to have a fixed point. 
In the case of the torsion-free metric connection (the Levi-Civita connection),the possible Riemannian holonomy groups 
are known to be $\mbox{SU}(n), \mbox{Sp}(n), \mbox{G}_2, \mbox{Spin}(7)$ \cite{Hit,Wang}. 
The Riemannian holonomy condition imposes strong restrictions on the geometry and leads 
to the consideration of Calabi-Yau manifolds, hyper-K\"ahler manifolds, parallel 
$\mbox{G}_2$-manifolds and parallel 
$\mbox{Spin}(7)$-manifolds. All of them are of great interest in mathematics 
(see \cite{J2} for detailed discussions) as well as in high-energy physics and
string theory \cite{Pol}. However, it seems that the geometry of these spaces is 
too restrictive for various problems in string theory \cite{MN,ST,GKMW}. 
One possible generalization of Calabi-Yau manifolds, hyper-K\"ahler manifolds, 
parallel $\mbox{G}_2$-manifolds and parallel $\mbox{Spin}(7)$-manifolds are manifolds equipped with linear metric connections having skew-symmetric torsion and holonomy contained in $\mbox{SU}(n), \mbox{Sp}(n), \mbox{G}_2, \mbox{Spin}(7)$. One remarkable fact is that the existence (in small dimensions) of a 
parallel spinor with respect to a metric connection $\nabla$ with 
skew-symmetric torsion determines the connection in a unique way if its 
holonomy group is a subgroup of $\mbox{SU}, \mbox{Sp}, \mbox{G}_2$, provided that some additional differential conditions on the structure are fulfilled  \cite{Stro,FI}, and always in dimension $8$ for a subgroup of the group $\mbox{Spin}(7)$ \cite{I}. The case of $16$-dimensional Riemannian manifolds with $\mbox{Spin}(9)$-structure was investigated in \cite{Fri}, homogeneous models are discussed in \cite{Agr}.
The existence of $\nabla$-parallel spinors in the dimensions $4$, $5$, $6$, $7$, $8$ is studied 
in \cite{Stro,DI,IP,FI,FI1,I}. In dimension $7$, the first consequence is that 
the manifold should be a $\mbox{G}_2$-manifold with an integrable 
$\mbox{G}_2$-structure \cite{FI}, i.\,e.,\ the structure group could be 
reduced to the 
group $\mbox{G}_2$ and the corresponding  $3$-form $\omega^3$ should obey 
$d*\omega^3=\theta\wedge *\omega^3$ for some special $1$-form $\theta$. 
In this paper we study solutions to the Killing spinor equations (\ref{ks}) 
in dimension $7$ and the geometry of integrable $\mbox{G}_2$-manifolds. We find a formula for  
the Riemannian scalar curvature in terms of the fundamental $3$-form. Our first main result 
is the following
\begin{thm}\label{th2}
%---------------------
Let $(M,g,\omega^3)$ be an integrable $\mbox{G}_2$-manifold with the 
fundamental $3$-form $\omega^3$. The Riemannian scalar curvature $\mathrm{Scal}^g$ is given in terms of the 
fundamental 3-form $\omega^3$ by
\begin{eqnarray}\label{sc1}
\mathrm{Scal}^g\, =\, \frac{1}{18}(d\omega^3,*\omega^3)^2 + 2||\theta||^2 -
\frac{1}{12}||T||^2+3\delta\theta \, , 
\end{eqnarray}
where $\theta$ and $T$ are the Lee form and the torsion of the unique 
$\mbox{G}_2$-connection given by 
\beq\label{eq0}
T\,=\,-*d\omega^3 + \frac{1}{6}(d\omega^3,*\omega^3)\cdot\omega^3 +  *(\theta\wedge\omega^3), \quad
\theta\,=\,-\frac{1}{3} *(*d\omega^3\wedge\omega^3)\, = \, \frac{1}{3}*
(\delta\omega^3\wedge *\omega^3)\,.
\eeq
\end{thm}
\noindent
We remark that the torsion form $T$ was been computed in \cite{FI}.
Returning to the Killing spinor equations (\ref{ks}), we 
present necessary and sufficient  conditions for a $\mbox{G}_2$-manifold to be 
a solution to  both of them. In fact we show that the dilation function 
arises from the Lee $1$-form. Finally, we give a formula for the Riemannian 
scalar curvature of any solution to  both Killing spinor equations
in dimension $7$. Our second main result is 
\begin{thm}\label{th3}
%---------------------
In dimension $7$ the following conditions are equivalent:
\begin{enumerate}
\item The Killing spinor equations (\ref{ks}) admit a solution with dilation 
$\Phi$;
\item There exists an integrable $\mbox{G}_2$-structure $(g,\omega^3)$ 
with closed Lee form, which is 
locally conformally equivalent to a cocalibrated $\mbox{G}_2$-structure of pure 
type $W_3$.
\end{enumerate}
More precisely, the structure is determined by the equations
\beq\label{5}
d*\omega^3\ =\ \theta\wedge *\omega^3, \qquad (d\omega^3,*\omega^3)\ = \ 0, 
\qquad \theta\ = \ -\, 2d\Phi
\eeq
and the NS $3$-form $H=T$ is given by
\beq\label{6}
T\ =\ -*d\omega^3  - 2*(d\Phi\wedge\omega^3).
\eeq
The Riemannian scalar curvature is determined by
\beq\label{lsc}
\mathrm{Scal}^g \ = \ 8\cdot ||d\Phi||^2 -\frac{1}{12}\cdot ||T||^2-6\cdot
\triangle\Phi \ ,
\eeq
where $\triangle \Phi=\delta d\Phi$ is the Laplacian. The solution has 
constant dilation if and only if the $\mbox{G}_2$-structure is cocalibrated of pure type $W_3$.
\end{thm}
\noindent
Our proof relies on the existence theorem for a $\mbox{G}_2$-connection with torsion, 
the Schr\"odinger-Lichnerowicz formula for the connection with torsion (both 
established in \cite{FI}) and the special properties of the Clifford action on the 
special parallel spinor.
%
%----------------------------------------------------------------------------------------
\section{General properties of $\mbox{G}_2$-structures}\noindent
%----------------------------------------------------------------------------------------
Let us consider $\mathbb{R}^7$ endowed with an orientation and its standard 
inner product. Denote an oriented orthonormal basis
by $e_1,\ldots,e_7$. We 
shall use the same notation for the dual basis. We denote the monomial 
$e_i\wedge e_j \wedge e_k$
by $e_{ijk}$. Consider the $3$-form $\omega^3$ on 
$\mathbb{R}^7$ given by
\begin{eqnarray}\label{1}
\omega^3 &=&e_{127} + e_{135} -
e_{146}-e_{236} - e_{245} +  
e_{347} + e_{567}.
\end{eqnarray}
The subgroup of $\mbox{SO}(7)$ that fixes $\omega^3$ is the exceptional Lie group 
$\mbox{G}_2$. It is a compact, simply-connected, simple Lie group of dimension $14$ 
\cite{Reichel}. The $3$-form $\omega^3$ corresponds to a real spinor and, 
therefore, $\mbox{G}_2$ is the isotropy group of a non-trivial real spinor.  A 
{\it $\mbox{G}_2$-structure} on a $7$-manifold $M^7$ is a reduction of the structure 
group of the tangent bundle to the exceptional group $\mbox{G}_2$. This can be described 
geometrically by a nowhere vanishing differential $3$-form $\omega^3$ on $M^7$, which can be 
locally written as (\ref{1}). The 3-form $\omega^3$ is called the {\it fundamental form} of 
the $\mbox{G}_2$-manifold $M^7$ (see \cite{Bo}) and it determines the metric completely. The 
action of $\mbox{G}_2$ on the tangent space 
gives an action of $\mbox{G}_2$ on $k$-forms and we obtain the following splitting 
\cite{FG,Br}:
\begin{displaymath}
\Lambda^1(M^7)\ = \ \Lambda^1_7, \quad \Lambda^2(M^7) \ = \ \Lambda^2_7\oplus 
\Lambda^2_{14}, \quad \Lambda^3(M^7)\ = \ \Lambda^3_1\oplus\Lambda^3_{7}\oplus
\Lambda^3_{27} \ ,
\end{displaymath}
where
\begin{displaymath}
\Lambda^2_7 \ = \ \{\alpha \in \Lambda^2(M^7)\ |\ *(\alpha\wedge\omega^3)=2\alpha\}, \quad
\Lambda^2_{14} \ = \ \{\alpha \in \Lambda^2(M^7)\|\ *(\alpha\wedge\omega^3)=
-\alpha\} \, ,
\end{displaymath}
\begin{displaymath}
\Lambda^3_7 \ = \ \{*(\beta\wedge\omega^3)\ |\ \beta \in \Lambda^1(M^7)\},\quad 
\Lambda^3_{27} \ = \ \{\gamma \in \Lambda^3(M^7)\ |\ \gamma\wedge\omega^3=0, \gamma\wedge *
\omega^3=0\}\,.
\end{displaymath} 
\noindent
Following \cite{Cab} we consider the $1$-form $\theta$ defined by
\beq\label{c2}
3\theta=- *(*d\omega^3\wedge\omega^3)=*(\delta\omega^3\wedge *\omega^3)\ .
\eeq
We shall call this $1$-form the {\it Lee form} associated with a given $\mbox{G}_2$-structure. 
If the Lee form vanishes, then we shall call the $\mbox{G}_2$-structure 
{\it balanced}. The classification of the different types of $\mbox{G}_2$-structures was worked out by Fernandez-Gray \cite{FG}, and Cabrera used the Lee 
form to characterize each of the $16$ classes. An {\it integrable} $\mbox{G}_2$-structure 
(or a structure of type $W_1\oplus W_3\oplus W_4$ ) is characterized by the differential 
equation
\begin{displaymath} 
d*\omega^3 \ = \ \theta\wedge * \omega^3 \, ,
\end{displaymath}
and a {\it cocalibrated} $\mbox{G}_2$-structure is defined by the condition
\begin{displaymath}
d * \omega^3 \ = \ 0 \, .
\end{displaymath}
\noindent
A {\it cocalibrated} $G_2$-{\it structure of pure type} $W_3$ is characterized by the two 
condtions $d* \omega^3 = 0, \ d \omega^3 \wedge \omega^3 = 0$. Then the 
following proposition follows immediately.
\begin{pro}\label{cor1}
If the Lee $1$-form is closed, then the $\mbox{G}_2$-structure is locally conformal to a 
balanced $\mbox{G}_2$-structure.
\end{pro}
\noindent
We shall call locally conformally parallel $\mbox{G}_2$-manifolds that are not globally 
conformally parallel {\it strict locally conformally parallel}.

\begin{exa}
Any $7$-dimensional oriented spin Riemannian manifold admits a certain 
$\mbox{G}_2$-structure, in general a non-parallel one (see for example \cite{LM}). 
The first known examples of complete parallel $\mbox{G}_2$-manifold were
constructed by Bryant and Salamon \cite{BS}, the first compact examples by 
Joyce \cite{J1,J2,J3}. There are many known examples of compact nearly parallel 
$\mbox{G}_2$-manifolds: $S^7$ \cite{FG}, $\mbox{SO}(5)/\mbox{SO}(3)$ \cite{BS,Sal}, 
the Aloff-Wallach spaces $N(g,l)=\mbox{SU}(3)/\mbox{U}(1)_{g,l}$ \cite{CMS}, any 
Einstein-Sasakian and any 3-Sasakian space in dimension $7$ \cite{FK,FKMS}. 
There are also some non-regular $3$-Sasakian manifolds (see \cite{6,7}). Moreover,
compact nearly parallel $\mbox{G}_2$-manifolds with large symmetry group are classified 
in \cite{FKMS}. Compact integrable nilmanifolds are constructed and studied in \cite{FU}.
Any minimal hypersurface $N$ in $\mathbb{R}^8$ admits a cocalibrated  
$\mbox{G}_2$-structure \cite{FG}. Moreover, the structure is parallel, nearly parallel, 
cocalibrated of pure type if and only if the hypersurface $N$ is totally geodesic, totally umbilic or minimal, respectively.
\end{exa}
%
%--------------------------------------------------------------------------------------------
\section{Conformal transformations of $\mbox{G}_2$-structures}\noindent
%--------------------------------------------------------------------------------------------
We study the conformal transformation of $\mbox{G}_2$-structures (see \cite{FG}).
\begin{pro}\label{pro1}
%----------------------
Let $\bar{g}=e^{2f}\cdot g, \ \bar{\omega^3}=e^{3f}\cdot \omega^3$  
be a conformal change of a $\mbox{G}_2$-structure $(g,\omega^3)$ and denote by 
$\bar{\theta}, \theta$ the corresponding Lee forms, respectively. Then
\beq\label{2}
\bar{\theta}\,=\, \theta + 4df.
\eeq
\end{pro}
\begin{proof}
We have the relations  
\begin{displaymath}
\mathrm{vol}_{\bar{g}}\ = \ e^{7f}\cdot \mathrm{vol}_g, \quad 
d\bar{\omega^3} \ = \ e^{3f} \cdot (3df\wedge\omega^3 + d\omega^3) \, . 
\end{displaymath}
We calculate 
$$\bar{*}d\bar{\omega^3}=e^{4f}(*d\omega^3 + 3*(df\wedge\omega^3)), \quad
\bar{*}d\bar{\omega^3}\wedge\bar{\omega^3}=e^{7f}(*d\omega^3\wedge\omega^3 - 12*df),$$
where we used the general identity $*(\omega^3\wedge\gamma)\wedge\omega^3
=4*\gamma$, which is valid for any 1-form $\gamma$. Consequently, we obtain
$\bar{\theta}= -\frac{1}{3}\bar{*}(\bar{*}d\bar{\omega^3}\wedge\bar{\omega^3})=-\frac{1}{3}\left(
*(*d\omega^3\wedge\omega^3) -12 *^2df\right)= \theta +4df.$ 
\end{proof} 
\noindent
Proposition~\ref{pro1} allows us to find a distinguished $\mbox{G}_2$-structure on a compact 7-dimensional $\mbox{G}_2$-manifold. 
\begin{thm}\label{th1}
Let $(M^7,g,\omega^3)$ be a compact 7-dimensional $\mbox{G}_2$-manifold. Then there 
exists a unique (up to homothety) conformal $\mbox{G}_2$-structure $g_0=e^{2f}\cdot g,
 \omega^3_0=e^{3f}\cdot \omega^3$ such that 
the corresponding Lee form is coclosed, $\delta_0\theta_0=0$.
\end{thm}
\begin{proof}
We shall use  the Gauduchon Theorem for the existence of a distinguished 
metric on a compact,  hermitian or Weyl manifold \cite{G1,G2}. We shall use 
the expression of this theorem in terms of a Weyl structure (see \cite{tod}, 
Appendix 1). We consider the Weyl manifold $(M^7,g,\theta,\nabla^W)$ with the 
Weyl 1-form $\theta$, where $\nabla^W$ 
is a torsion-free linear connection on $M^7$ determined by the condition
$\nabla^Wg=\theta\otimes g$. Applying the Gauduchon Theorem we can find, in a 
unique way, a conformal metric $g_0$ such that the corresponding Weyl 1-form 
is coclosed with respect to $g_0$. The key point is that, by Proposition~\ref{pro1}, the Lee 
form transforms under conformal rescaling according to (\ref{2}), which is exactly the 
transformation of the Weyl $1$-form 
under conformal rescaling of the metric $\bar g=e^{4f}\cdot g$. 
Thus, there exists (up to homothety) a unique  conformal
$\mbox{G}_2$-structure $(g_0,\omega^3_0)$ with coclosed Lee form. 
\end{proof}
\noindent
We shall call the $\mbox{G}_2$-structure with coclosed Lee form 
{\it the Gauduchon $\mbox{G}_2$-structure}.
\begin{co}
%---------
Let $(M^7,g,\Phi)$ be a compact $\mbox{G}_2$-manifold and $(g,\Phi)$ be the Gauduchon 
structure. Then the following formula holds:
\begin{displaymath}
*\left(d\delta\omega^3\wedge *\omega^3\right) \ = \ ||\delta\omega^3||^2 \, .
\end{displaymath}
In particular, if the structure is integrable, then 
\begin{displaymath}
*\left(d\delta\omega^3\wedge *\omega^3\right) \ = \ 24||\theta||^2 \ .
\end{displaymath}
\end{co}
\begin{proof} 
%------------
Using (\ref{c2}), we calculate that
\begin{displaymath}
0 \ = \ 3 \cdot \delta\theta =  *d(\delta\omega^3\wedge *\omega^3) =  
*\left(d\delta\omega^3\wedge *\omega^3 - *d*\omega^3\wedge d*\omega^3\right)
=  *\left(d\delta\omega^3\wedge *\omega^3 -||\delta\omega^3||^2 \cdot \mathrm{vol}
\right) \, .
\end{displaymath}  
If the structure is integrable, then  $\ ||\delta\omega^3||^2=24||\theta||^2$. 
\end{proof}
\begin{co}\label{col}
On a compact $\mbox{G}_2$-manifold with closed Lee form whose Gauduchon $\mbox{G}_2$-structure
is not balanced, the first Betti number satisfies
$b_1(M)\ge 1$.
\end{co}
\noindent
For integrable $G_2$-manifolds one can define a suitable elliptic complex
as well as cohomolgy groups $\tilde{\mbox{H}}^i(M^7)$ (see \cite{FU}). The first cohomolgy group is given by
\begin{displaymath}
\tilde{\mbox{H}}^1(M^7) \ = \ \{ \alpha \in \Lambda^1(M^7) : d \alpha \wedge 
* \omega^3 = 0, \quad d*\alpha = 0 \} \ .
\end{displaymath}
\begin{co}
On a compact integrable manifold which is not globally conformally balanced, one has
$\tilde b_1\ge 1$.
\end{co}
\begin{proof} By the condition of the theorem the Gauduchon structure has a 
non-identically zero
Lee form. Then $0=\delta\omega^3=*(d\theta\wedge *\omega^3)$, since the structure is integrable.
Adding the condition $\delta\theta=0$, we obtain $\tilde b_1\ge 1$. 
\end{proof} 

\section{Connections with torsion, parallel spinors and Riemannian scalar curvature}
%------------------------------------------------------------------------------------
\noindent
The Ricci tensor of an integrable $\mbox{G}_2$-manifold was expressed in principle by the
structure form $\omega^3$ in the paper \cite{FI}. Here we intend to find 
an explicit formula for the Riemannian scalar curvature.
Using the unique connection with skew-symmetric torsion
preserving the given integrable $\mbox{G}_2$-structure found in \cite{FI}, we apply the 
Schr\"odinger-Lichnerowicz formula for the Dirac operator of a metric connection with 
totally skew-symmetric torsion (see \cite{FI}) in order to derive the formula for 
the scalar curvature. First, let us summarize the mentioned results from \cite{FI}.
\begin{thm}\label{tA}(see \cite{FI})
%-----------------------------------
Let $(M^7,g,\omega^3)$ be a $\mbox{G}_2$-manifold. Then the following conditions 
are equivalent:
\begin{enumerate}
\item The $\mbox{G}_2$-structure is integrable, i.e., $d*\omega^3=\theta\wedge *\omega^3 $;
\item There exists a unique linear connection $\nabla$ preserving the $\mbox{G}_2$-structure 
with totally skew-symmetric torsion $T$ given by
\beq\label{3}
T \ = \ -\, *d\omega^3 + \frac{1}{6}(d\omega^3,*\omega^3)\cdot\omega^3 +  
*(\theta\wedge\omega^3)\ .
\eeq
\end{enumerate}
Furthermore, for any integrable $\mbox{G}_2$-structure, the projections
 $\pi^4_1(d\omega^3), \ \pi^4_7(d\omega^3)$ of $d\omega^3$ onto 
$\Lambda^4_1$ and $\Lambda^4_7$, respectively, are given by
\begin{displaymath}
\pi^4_1(d\omega^3) \ = \ \frac{1}{7} \cdot (d\omega^3,*\omega^3)*\omega^3, \quad \pi^4_7(d\omega^3) \ = \ \frac{3}{4}\cdot \theta\wedge\omega^3 \ ,
\end{displaymath}
there exists a $\nabla$-parallel spinor $\Psi_0$ corresponding to the fundamental form 
$\omega^3$ and the Clifford action of the torsion $3$-form on it is 
\beq\label{b2}
T\cdot \Psi_0 \ = \ \frac{7}{6}\cdot \lambda \cdot \Psi_0 - \theta \cdot \Psi_0, \quad 
\lambda \ = \ - \, \frac{1}{7}\cdot(d\omega^3,*\omega^3) \ .
\eeq
\end{thm}
\noindent
Keeping in mind Proposition~\ref{pro1}, we obtain
\begin{co}\label{co2}
The Lee form of an integrable $\mbox{G}_2$-structure is given by 
 $\ *(\omega^3\wedge T)= -\theta $.
\end{co}
\begin{co}
The torsion $3$-form $T$ of $\nabla$ changes by a conformal transformation 
$(g_o=e^{2f}\cdot g, \omega^3_o=e^{3f}\cdot\omega^3)$ of the $\mbox{G}_2$-structure by 
\begin{displaymath}
T_o=e^{4f}\cdot(T + *df\wedge\omega^3).
\end{displaymath}
\end{co} 
\noindent
Let $\mathrm{D}$ and $\mathrm{Scal}$ be the Dirac operator  and the 
scalar curvature of the $\mbox{G}_2$-connection defined as usually by
\begin{displaymath}
\mathrm{D} \ = \ \sum_{i=1}^7 e_i \cdot \nabla_{e_i},\quad
\mathrm{Scal} \  = \ \sum_{i,j=1}^7R^{\nabla}(e_i,e_j,e_j,e_i) \ .
\end{displaymath}
The scalar curvature $\mathrm{Scal}^g$ of the metric is given by (see \cite{IP,FI})
\beq\label{c5}
\mathrm{Scal}^g \ = \ \mathrm{Scal} + \frac{1}{4}||T||^2 \ .
\eeq
The $4$-form $\sigma^T$ defined by the formula
\begin{displaymath}
\sigma^T \ = \ \frac{1}{2}\sum_{i=0}^7(e_i \haken T)\wedge (e_i \haken T)  
\end{displaymath}
plays an important role in the integrability conditions for $\nabla$-parallel spinors.
\begin{thm}\label{tB}(see \cite{FI}) 
%-----------------------------------
Let $\Psi$ be a parallel spinor with respect to a metric connection $\nabla$ with totally skew-symmetric torsion $T$ on a Riemannian spin manifold $M^n$. 
Then the following formulas hold
\begin{displaymath}\label{c3}
3 \cdot dT\cdot \Psi -2 \cdot \sigma^T\cdot \Psi + \mathrm{Scal}\cdot \Psi \ = \ 0, \quad 
\mathrm{D}(T\cdot\Psi) \ = \ dT\cdot\Psi +\delta T\cdot\Psi -2\cdot\sigma^T\cdot\Psi \, .
\end{displaymath}
\end{thm}

\noindent
\textit{Proof of Theorem $\ref{th2}$.}
Let $\Psi_0$ be the $\nabla$-parallel spinor corresponding to the fundamental 
3-form $\omega^3$. Then the Riemannian Dirac operator $\mathrm{D}^g$ and the Levi-Civita connection $\nabla^g$ act on $\Psi_0$ by the rule
\beq\label{dir1}
\nabla^g_X\Psi_0 \  = \ -\, \frac{1}{4}(X \haken T)\cdot \Psi_0, \quad 
\mathrm{D}^g\Psi_0 \  = \  -\, \frac{3}{4}\cdot T\cdot \Psi_0 \ = \
-\, \frac{7}{8}\cdot \lambda \cdot \Psi_0 +\frac{3}{4}\cdot \theta\cdot \Psi_0 \, ,
\eeq
where we used Theorem~\ref{tA}. We are going to apply the well known 
Schr\"odinger-Lichnerowicz formula \cite{Li,Sro}
$$(\mathrm{D}^g)^2 \ = \ \triangle^g  + \frac{1}{4}\cdot \mathrm{Scal}^g, 
\quad \triangle^g \ = \  -\, \sum_{i=1}^n\left( \nabla^g_{e_i}\nabla^g_{e_i} -\nabla^g_{\nabla_{e_i}e_i}\right)
$$
to the  $\nabla$-parallel spinor field $\Psi_0$. The formula (\ref{dir1}) yields that
\begin{eqnarray}\label{d2}
& &(\mathrm{D}^g)^2\Psi_0 \ = \ -\, \frac{7}{8}\cdot \mathrm{D}^g\big(\lambda
\cdot \Psi_0\big) + \frac{3}{4}\cdot \mathrm{D}^g\big(\theta\cdot \Psi_0\big)\ =\ \\ \nonumber
& & = \ \left(\frac{49}{64}\cdot \lambda^2+\frac{9}{16}\cdot ||\theta||^2+\frac{3}{4}\cdot 
\delta\theta\right)\cdot \Psi_0 -\frac{7}{8}\cdot d\lambda\cdot \Psi_0 + \frac{3}{4}\cdot 
d\theta\cdot \Psi_0 + \frac{3}{8}\cdot ({\theta} \haken T)\cdot \Psi_0 \, ,
\end{eqnarray}
where we used the general identity 
$\mathrm{D}^g\circ\theta+\theta\circ\mathrm{D}^g=d\theta+\delta\theta-2\nabla_{\theta}$. We 
compute the Laplacian $\triangle^g$. Fix a normal coordinate system at a point $p \in M^n$ 
such that $(\nabla_{ei}e_i)_p=0$, use (\ref{dir1}) as well as the properties of the 
Clifford multiplication. Then one obtains the following formula \cite{I}:
\begin{eqnarray}\label{d3}
& &\triangle^g\Psi_0 \ = \ \frac{1}{4}\cdot \sum_{i=1}^n\left(\nabla_{e_i}i_{e_i}T)\cdot \Psi_0 -\frac{1}{16}\cdot (e_i \haken T)\cdot 
(e_i \haken T)\cdot \Psi_0 \right)\\
\nonumber
& & = \ -\, \frac{1}{4}\cdot \delta T\cdot 
\Psi_0 -\frac{1}{16}\cdot \left(2\sigma^T-\frac{1}{2}\cdot ||T||^2 \right)\cdot \Psi_0 \, .
\end{eqnarray}
Substituting (\ref{d2}) and (\ref{d3}) into the SL-formula, multiplying the obtained result 
by $\Psi_0$ and taking the real part, we arrive at 
\beq\label{d4}
\left(\frac{49}{64}\cdot \lambda^2+\frac{9}{16}\cdot ||\theta||^2+\frac{3}{4}
\cdot \delta\theta\right)\cdot ||\Psi_0||^2 \ = \ 
\left(\frac{1}{32}^2+\frac{1}{4}\cdot \mathrm{Scal}^g\right)\cdot ||\Psi_0||^2
-\frac{1}{8}\cdot(\sigma^T\cdot \Psi_0,\Psi_0) \, .
\eeq
On the other hand, using (\ref{b2}), we obtain  
$$\mathrm{D}(T\cdot \Psi_0) \ = \ \mathrm{D}(\frac{7}{6}\cdot \lambda \cdot 
\Psi_0 - \theta\cdot \Psi_0) \ = \ \sum_{i=1}^n e_i \cdot \nabla_{e_i}
(\frac{7}{6}\cdot \lambda \cdot \Psi_0  - \theta \cdot \Psi_0) \ = \ 
\frac{7}{6}\cdot d\lambda \cdot \Psi_0 - \left(d^{\nabla}\theta + \delta\theta\right)
\cdot \Psi_0,$$ 
where $d^{\nabla}$ is the exterior
derivative with respect to the $\mbox{G}_2$-connection $\nabla$. Now,  
Theorem (\ref{c3}) gives $\frac{7}{6}d\lambda \cdot \Psi_0 -
d^{\nabla}\theta \cdot \Psi_0 - \delta\theta \cdot \Psi_0 = dT \cdot \Psi_0 
-2\sigma^T \cdot \Psi_0 + \delta T \cdot \Psi_0$.
Multiplying the latter equality by $\Psi_0$ and taking the real part, we obtain
$-\delta\theta \cdot ||\Psi_0||^2=(dT\cdot\Psi_0,\Psi_0) -(2\sigma^T\cdot
\Psi_0,\Psi_0)$. Consequently, Theorem (\ref{c3}) and (\ref{c5}) imply
\beq\label{d5}
\left(-3\cdot \delta\theta-\frac{1}{4}\cdot||T||^2 + \mathrm{Scal}^g\right)
\cdot ||\Psi_0||^2 + 4\cdot(\sigma^T \cdot \Psi_0,\Psi_0) \ = \ 0 \, .
\eeq
Finally, (\ref{d4}) and (\ref{d5}) imply (\ref{sc1}) and the proof of 
Theorem~\ref{th2} is complete. \qed
\begin{co}
On a cocalibrated $\mbox{G}_2$-manifold of pure type the Riemannian scalar 
curvature is given by 
\begin{displaymath}
\mathrm{Scal}^g\ = \ - \, \frac{1}{12}\cdot ||d\omega^3||^2 \ .
\end{displaymath}
\end{co}
\begin{proof}
%------------
In the case of a cocalibrated $\mbox{G}_2$-structure of pure type, the 
torsion 3-form 
$T=-*d\omega^3$. The claim follows from Theorem~\ref{th2}.
\end{proof} 
\noindent
Using the results in \cite{FG} we derive immediately the following formula, 
which is essentially the reformulated Gauss equation.
\begin{co}\label{min}
Let $M^7$ be a hypersurface in $\mathbb{R}^8$ the with second fundamental form 
$\mathrm{S}$ and mean curvature $\mathrm{H}$. Then the Riemannian scalar curvature on 
$M^7$ is given by the formula
\begin{eqnarray}\label{sur}
\mathrm{Scal}^g \ = \ \frac{49}{18}\cdot||\mathrm{H}||^2 -\frac{1}{12}\cdot ||
\mathrm{S}_0||^2, 
\end{eqnarray}
where $\mathrm{S}_0$ is the image of the traceless part of the second fundamental form via the isomorphism $S^2_0(R^7)\rightarrow \Lambda^3_{27}$. In 
particular, if $M$ is a minimal hypersurface, then
\begin{eqnarray}\label{min1}
\mathrm{Scal}^g \ = \  -\frac{1}{12}||S_0||^2 \ \le \ 0 \ .
\end{eqnarray}
\end{co}

\begin{thm}\label{mt}
Let $M^7$ be a compact, connected spin $7$-manifold with a fixed orientation. 
If it admits a strictly locally conformally parallel $\mbox{G}_2$-structure, then:
\begin{enumerate}
\item $M$ admits a Riemannian metric $g_Y$ with strictly positive constant 
scalar curvature,
\item the first Betti number is at least one, $b_1(M)\ge 1$.
\end{enumerate}
\end{thm}
\begin{proof} 
%------------
We have $||T||^2=\frac{3}{2}||\theta||^2$ since the structure is locally 
conformally parallel. Then, Theorem~\ref{th2} leads to the formula
\beq\label{as}
\mathrm{Scal}^g \ = \ \frac{15}{8}\cdot ||\theta||^2 +3\cdot \delta\theta \, .
\eeq
According to the solution of the Yamabe conjecture \cite{RS} there is a metric 
$g_Y=e^{2f}\cdot g$ in the conformal class of $g$ with constant scalar curvature.
Consider the locally conformally parallel $\mbox{G}_2$-structure 
$(g_Y=e^{2f}\cdot g, \omega^3_Y=e^{3f}\cdot \omega^3)$. The 
equality (\ref{as}) also holds for the structure $(g_Y,\omega^3_Y)$ and  an integration 
over $M$ gives
$$ \mathrm{Scal}^{g_Y}\cdot \mathrm{vol}(g_Y) \ = \ \frac{11}{6}\int_M||\theta||^2\,d\mathrm{vol} \ > \ 0 \ ,$$
since the structure is strictly locally conformally parallel. The second assertion is a consequence of Corollary~\ref{col}.
\end{proof}

\section{Solutions to the Killing spinor equations in dimension 7}\noindent
%---------------------------------------------------------------------------
We consider the Killing spinor equations (\ref{ks}) in dimension $7$. 
The existence of a $\nabla$-parallel spinor is equivalent to  the existence of 
a $\nabla$-parallel integrable $\mbox{G}_2$-structure and the $3$-form 
field strength $H=T$ is given by (\ref{3}). We now investigate the second Killing spinor 
equation (\ref{ks}).\\

\noindent
\textit{Proof of Theorem $\ref{th3}$}.
Let $\Psi$ be an arbitrary $\nabla$-parallel such that $(d \Phi - T) \cdot \Psi = 0$. The 
spinor field $\Psi$ defines a second $\mbox{G}_2$-structure $\omega^3_0$ such that 
$\Psi = \Psi_0$ is the canonical spinor field. Since
the connection preserves the spinor field $\Psi$, it preserves the $\mbox{G}_2$-structure
 $\omega^3_0$, too. On the other hand, the connection preserving $\omega^3_0$ is unique. 
Consequently, the torsion $T_0$ coincides with the torsion form
$T$ and for the $\mbox{G}_2$-structure $\omega^3_0$ we have
\begin{displaymath}
\nabla \Psi_0 \ = \ 0, \quad (d \Phi - T_0) \cdot \Psi_0 \ = \ 0 \, .
\end{displaymath}
The Clifford action $T_0 \cdot \Psi_0$ depends only on the 
$(\Lambda^3_1\oplus\Lambda^3_7)$-part of $T_0$. Using (\ref{3}) and the algebraic formulas 
\begin{displaymath}
*(\gamma\wedge\omega^3_0) \cdot \Psi_0 \ = \ - \gamma \haken (*\omega^3_0) \cdot \Psi_0 \ = \ - \, 4\cdot \gamma \cdot \Psi_0 ,\quad 
\omega^3_0 \cdot \Psi_0 \ = \ - \, 7 \cdot \Psi_0 
\end{displaymath}
we calculate
\beq\label{a3}
T_0 \cdot \Psi_0 \  = \  - \, \theta \cdot \Psi_0 - \frac{1}{6}\cdot 
(d\omega^3_0,*\omega^3_0) \cdot \Psi_0 \ .
\eeq
Comparing with the second Killing spinor equation (\ref{ks}) we find 
$2 \cdot d\Phi=-\beta, \ (d\omega^3_0,*\omega^3_0)=0$ which completes the proof. \qed \\

\noindent
As a corollary we obtain the result from \cite{GKMW}, which states that any solution to 
both equations (\ref{ks}) has necessarily the NS three form $H=T$ given 
by (\ref{6}). A more precise analysis using Proposition~\ref{pro1} and 
Theorem~\ref{th2} of the explicit solutions constructed in \cite{GKMW} shows 
that these solutions are conformally equivalent to a cocalibrated structure of pure type. In 
other words, the multiplication of the $\mbox{G}_2$-structures $(g^{\pm},\omega^{3\pm})$ by 
$(e^{\Phi}\cdot g^{\pm},e^{(3/2)\Phi}\cdot \omega^{3\pm})$ 
is a new example of a cocalibrated $\mbox{G}_2$-structure of pure type $W_3$, 
and it is a solution to the Killing spinor equations with constant dilation. The same 
conclusions are valid for the solutions constructed in \cite{AGK,ST,MN}.\\

\noindent
Theorem~\ref{th3} allows us to construct a lot of compact solutions to the Killing spinor 
equations. If the dilation is a globally defined function, 
then any solution is globally conformally equivalent to a cocalibrated $\mbox{G}_2$-structure of pure type. For example, any conformal transformation of a compact 
$7$-dimensional manifold with a Riemannian holonomy group $\mbox{G}_2$ 
constructed by Joyce \cite{J1,J2} is a solution with non-constant dilation. 
Another source of solutions are conformal transformations of the cocalibrated 
$\mbox{G}_2$-structures of pure type $W_3$ induced on any minimal  hypersurface 
in $\mathbb{R}^8$. Summarizing, we obtain:
\begin{co}
Any solution $(M^7,g,\omega^3)$ to the Killing spinor equations (\ref{ks}) 
in dimension $7$ with non-constant globally defined dilation function $\Phi$ 
comes from 
a solution with constant dilation by a conformal transformation 
$(g=e^{\Phi} \cdot g_0,\omega^3=e^{(3/2)\Phi} \cdot \omega^3_0)$,
where $(g_0,\omega^3_0)$ is a cocalibrated $\mbox{G}_2$-structure of pure type $W_3$.
\end{co}


\begin{thebibliography}{99}
%\bibliographystyle{amsalpha}
\bibitem {AGK} B. Acharya, J. Gauntlett, N. Kim, {\em Fivebranes wrapped on associative three-cycles}, to appear in Phys. Rev. D, hep-th/0011190.
\bibitem {Agr} I. Agricola, {\em Connections on naturally reductive spaces,
their Dirac operator and homogeneous models in string theory}, to appear.
\bibitem {Bo} E. Bonan, {\em Sur le vari\'et\'es riemanniennes a groupe 
d'holonomie $\mbox{G}_2$ ou $\mbox{Spin}(7)$}, C. R. Acad. Sci. Paris {\bf 262} (1966), 
127-129.
\bibitem{6} C. Boyer, K. Galicki, B. Mann, {\em Quaternionic reduction and 
Einstein manifolds}, Comm. Anal. Geom., {\bf 1} (1993), 1-51.
\bibitem{7}  C. Boyer, K. Galicki, B. Mann, {\em The geometry and topology of 
3-Sasakian manifolds}, J. reine ang. Math. {\bf 455} (1994), 183-220.
\bibitem{Br} R. Bryant, {\em Metrics with exceptional holonomy}, Ann. Math. 
{\bf 126} (1987), 525-576.
\bibitem{BS} R. Bryant, S. Salamon, {\em On the construction of some complete 
metrics with exceptional holonomy}, Duke Math. J. {\bf 58} (1989), 829-850.
\bibitem {Cab} F. Cabrera, {\em On Riemannian manifolds with $\mbox{G}_2$-structure}, Bolletino U.M.I. (7) 10-A (1996), 98-112.
\bibitem{CMS} F. Cabrera, M. Monar, A. Swann, {\em Classification of $\mbox{G}_2$-structures}, J. London Math. Soc. {\bf 53} (1996), 407-416.
\bibitem{DI} P. Dalakov, S. Ivanov, {\em Harmonic spinors of Dirac operator of 
connection with torsion in dimension 4}, Class. Quantum Gravity {\bf 18} (2001), 253-265.
\bibitem {FG} M. Fernandez, A. Gray, {\em Riemannian manifolds with structure 
group $\mbox{G}_2$}, Ann. Mat. Pura Appl. {\bf 32} (1982), 19-45.
\bibitem {FU} M. Fernandez, L. Ugarte, {\em Dolbeault cohomology for $\mbox{G}_2$-manifolds}, Geom. Dedicata, {\bf 70} (1998), 57-86.
\bibitem{Fri} Th. Friedrich, {\em $\mbox{Spin}(9)$-structures and connections
with totally skew-symmetric torsion}, to appear.
\bibitem{FK} Th. Friedrich, I. Kath, {\em Compact 7-dimensional manifolds with Killing spinors}, Comm. Math. Phys. {\bf 133} (1990), 543-561.
\bibitem{FKMS} Th. Friedrich, I. Kath, A. Moroianu, U. Semmelmann, {\em On nearly parallel $\mbox{G}_2$-structures}, J. Geom. Phys. {\bf 23} (1997), 256-286.
\bibitem{FI} Th. Friedrich, S. Ivanov, {\em Parallel spinors and connections with skew symmetric torsion in string theory}, math.DG/0102142.
\bibitem{FI1} Th. Friedrich, S. Ivanov, {\em Almost contact manifolds and type II string equations}, math.DG/0111131.
\bibitem{GS} K. Galicki, S. Salamon, {\em On Betti numbers of 3-Sasakian manifolds}, Geom. Dedicata {\bf 63} (1996), 45-68.
\bibitem {G1} P. Gauduchon, {\em La 1-forme de torsion d'une vari\'et\'e
hermitienne compacte}, Math. Ann. {\bf 267}, (1984), 495-518.
\bibitem {G2} P. Gauduchon, {\em Structures de
 Weyl-Einstein,  espaces de twisteurs et
 vari\'{e}t\'{e}s de type $S^{1}\times S^{3}$}, J. reine ang. Math. 
{\bf 469}
 (1995), 1-50.
\bibitem {GKMW} J. Gauntlett, N. Kim, D. Martelli, D. Waldram, {\em Fivebranes wrapped on SLAG three-cycles and related geometry}, hep-th/0110034.
\bibitem {Gr} A. Gray, {\em Vector cross product on manifolds}, Trans. Am. 
Math. Soc. {\bf 141} (1969), 463-504, Correction {\bf 148} (1970), 625.
\bibitem{I} S. Ivanov, {\em Connections with torsion, parallel spinors and geometry of $\mbox{Spin}(7)$-manifolds}, math.DG/0111216
\bibitem{IP} S. Ivanov, G. Papadopoulos, {\em Vanishing theorems and string background}, Class. Quant. Grav. {\bf 18} (2001), 1089-1110.
\bibitem{IP1} S. Ivanov, G. Papadopoulos, {\em A no-go theorem for string warped compactification}, Phys. Lett. B {\bf 497} (2001) 309-316.
\bibitem{J1} D. Joyce, {\em Compact Riemannian 7-manifolds with holonomy $\mbox{G}_2$. I}, J.Diff. Geom. {\bf 43} (1996), 291-328.
\bibitem{J2} D. Joyce, {\em Compact Riemannian 7-manifolds with holonomy $\mbox{G}_2$. II}, J.Diff. Geom. {\bf 43} (1996), 329-375.
\bibitem{J3} D. Joyce, {\em Compact Riemannian manifolds with special holonomy}, Oxford University Press, 2000.
\bibitem {LM} B. Lawson, M.-L. Michelsohn, {\em Spin Geometry}, Princeton University Press, 1989.
\bibitem{Hit} N. Hitchin, {\em Harmonic spinors}, Adv. in Math. {\bf 14} (1974), 1-55.
\bibitem{Li} A. Lichnerowicz, {\it Spineurs harmoniques}, C. R. Acad. Sci. Paris, {\bf 257} (1963), 7-9.
\bibitem {MN} J. Maldacena, H. Nastase, {\em The supergravity dual of a theory with dynamical supersymmetry breaking}, JHEP 0109, {\bf 024} (2001), hep-th/0105049.
\bibitem {Pol} J. Polchinski, {\em String Theory vol.II, Superstring Theory and Beyond}, Cambridge Monographs on Mathematical Physics, Cambridge, University Press, 1998.
\bibitem{Reichel} W. Reichel, {\em \"Uber die Trilinearen alternierenden
Formen in $6$ und $7$ Variablen}, Dissertation Univ. Greifswald 1907.
\bibitem{Sal} S. Salamon, {\em Riemannian geometry and holonomy groups}, Pitman Res. Notes Math. Ser., 201 (1989).
\bibitem{RS} R. Schoen, {\em Conformal deformations of Riemannian metrics to constant scalar curvature}, J. Diff. Geom. {\bf 20} (1984), 479-495.
\bibitem {ST} M. Schvelinger, T. Tran, {\em Supergravity duals of gauge field theories from \mbox{SU}(2)xU(1)gauge supergravity in five dimensions}, JHEP 0106, {\bf 025} (2001), hep/th0105019.
\bibitem{Sro} E. Schr\"odinger, {\it Diracsches Elektron im Schwerfeld I}, Sitzungberichte der Preussischen Akademie der Wissenschaften Phys.-Math. Klasse 1932, Verlag der Akademie der Wissenschaften, Berlin 1932, 436-460.
\bibitem{Stro} A. Strominger, {\it Superstrings with torsion}, Nucl. Physics {\bf B 274} (1986), 254-284.
\bibitem {tod} K.P. Tod, {\em Compact 3-dimensional Einstein-Weyl
structures}, J. London Math. Soc. {\bf 45} (1992), 341-351.
\bibitem{Wang} M. Wang, {\em Parallel spinors and parallel forms}, Ann. Glob. Anal. Geom. {\bf 7} (1989), 59-68.

\end{thebibliography}
\end{document}